\documentclass[11pt]{article}
\usepackage{graphicx}
\usepackage{balance} 
\usepackage{lipsum} 
\usepackage[sc]{mathpazo} 
\usepackage[T1]{fontenc} 
\usepackage{microtype} 
\usepackage[hmarginratio=1:1,top=30mm,columnsep=15pt]{geometry} 
\usepackage{multicol} 
\usepackage[hang, small,labelfont=bf,up,textfont=it,up]{caption} 
\usepackage{booktabs} 
\usepackage{float} 
\usepackage{hyperref} 
\usepackage{lettrine} 
\usepackage{paralist} 
\usepackage{bbm}
\usepackage{float}
\usepackage{amsmath,amssymb}  
\usepackage{color} 
\usepackage[square]{natbib}
\usepackage{abstract} 

\usepackage{titlesec} 
\titleformat{\section}[block]{\large\scshape\centering}{\thesection.}{1em}{} 
\titleformat{\subsection}[block]{\large}{\thesubsection.}{1em}{} 

\usepackage{fancyhdr} 
\newtheorem{theo}{Theorem}

\newtheorem{app}{Application}

\newtheorem{coro}{Corollary}

\newtheorem{defi}{Definition}
\newtheorem{exa}{Example}

\newtheorem{lem}{Lemma}
\newtheorem{nota}{Notation}

\newtheorem{prop}{Proposition}

\newtheorem{rem}{Remark}

\newcommand{\RR}{\mathbb{R}}

\newcommand{\PP}{\mathbb{P}}

\newcommand{\ZZ}{\mathbb{Z}}
\newcommand{\NN}{\mathbb{N}}

\newcommand{\pp}{\mathbb{P}}
\newcommand{\EE}{\mathbb{E}}

\newcommand{\mc}{\mathcal}

\newcommand{\mb}{\mathbb}
\newcommand{\lfled}{\longrightarrow}

\newcommand{\­}{\neq}

\newcommand{\²}{\leqslant}
\newcommand{\³}{\geqslant}
\newcommand{\Cov}{\mbox{Cov}}
\newcommand{\Lip}{\mbox{Lip}}
\newcommand{\Var}{\mbox{Var}}

\newcommand{\e}{\mathbf{e}}
\newcommand{\h}{\widehat}
\newcommand{\1}{\mb{1}}

\title{\vspace{-15mm}\fontsize{25pt}{10pt}\selectfont\textbf{Empirical central limit theorem  for cluster functionals without mixing.}\thanks{This research has been conducted as part of the project Labex MME-DII (ANR11-LBX-0023-01).}} 
\author{\large \textsc{Paul Doukhan\footnote{Institut Universitaire de France (IUF) and ${\ddag}$.}  and Jos\'e-Gregorio G\'omez\footnote{D\'epartement de Math\'ematiques, Universit\'e de Cergy - Pontoise. 95000 Cergy - Pontoise, France.}}}
\date{}

\begin{document}
\maketitle 
\footnotetext{E-mail addresses: \href{mailto:doukhan@u-cergy.fr}{doukhan@u-cergy.fr} and \href{mailto:jose.gomez@u-cergy.fr}{jose.gomez@u-cergy.fr}}

\thispagestyle{fancy} 


\begin{abstract} We prove central limit theorems (CLT) for empirical processes of extreme values cluster functionals as in Drees and Rootz\'en (2010). We use coupling properties enlightened for Dedecker \& Prieur's $\tau-$dependence coefficients in order to improve the conditions of dependence and continue to obtain these CLT. The assumptions are precisely set for particular processes and cluster functionals of interest. The number of excesses provides a complete example of a cluster functional for a simple non-mixing model (AR(1)-process)  for which ours results are definitely needed. We also give the expression explicit the covariance structure of  limit Gaussian process.

Also we include in this paper some results of Drees (2011) for the extremal index and some simulations for this index to demonstrate the accuracy of this technique.
\begin{quotation}
\end{quotation}
\end{abstract}



\section{Introduction}
  
 \citet{Drees2010}'s scheme prove limit theorems for empirical process of cluster functionals (EPCF). 
 In statistics, \citet{Gomez2015} proves a CLT in finite dimensional convergence (fidis) under weaker conditions, and considers an example for which a functional CLT is not necessary.
 
 We extend the result under $\tau$-weak dependence. The classical example of a non-mixing autoregressive model demonstrates the importance of this functional extension. Moreover the estimation of the extremal index provides us with a suitable example of application of the functional CLT. 

For a real-valued random process $(X_i)_{i\in\NN}$,  a typical example of a EPCF is the tail empirical process:
\begin{align}\label{5}
T_n(x)=\dfrac{1}{\sqrt{n v_n}}\sum_{i=1}^n\big(\mb{1}\{X_i>a_n x +u_n\}-\pp(X_1>a_n x +u_n)\big),\qquad x\³0,
\end{align}
where $v_n=\PP\{ X_{n,1}\­0\}$ is decreasing to zero, $(u_n)_{n\in\NN}$ is a non-decreasing sequence of thresholds and $(a_n)_{n\in\NN}$ is a sequence of positive constants.  
This process has beed considered by Drees and Rootz\'en under suitable dependence conditions (in particular under $\alpha$ and $\beta$-mixing conditions), where they prove its uniform convergence to a Gaussian process $T$ under additional conditions. For example, they prove the convergence of this tail empirical process for  the cases of $k$-dependent sequences or stable AR(1)-processes \citep{Rootzen1995},  ARCH(1)-processes \citep{Drees2002,Drees2003} and some applications for solutions of stochastic difference equations \citep{Drees2000,Drees2002,Drees2003}. Finally, they use the  cluster functionals setting in \citep{Yun2000} and \citep{Segers2003} to generalize such empirical processes under $\beta$-mixing in \citep{Drees2010}.
 
 Unfortunately, the mixing processes family is very restrictive. This can be noted with the following AR(1)-process, solution of the recursion: 
\begin{align}\label{AR}
X_k=\dfrac{1}{b}\big(X_{k-1}+\xi_k\big),\quad k\in\ZZ,
\end{align}
 where $b\³2$ is an integer and $(\xi_k)_{k\in\NN}$ are independent and uniformly distributed random variables on the set $U(b):=\{0,1,\ldots, b-1\}$ which is not even $\alpha-$mixing, as this  is shown in \citep{Andrews}  for $b=2$ and  in \citep{Ango} for $b>2$. Thus, {\it the results in \citep{Drees2010} can not be used here!}. However, such process  (\ref{AR}) is $\tau$-weakly dependent as is shown in \citep{Dedecker2004a}. The same situation happens in a general way for the causal Bernoulli shifts, Markovian models, etc. This is thus useful to improve on the CLT for empirical processes of extreme cluster functionals proposed by \citet{Drees2010} for  more general classes of weakly dependent processes. 
 
In order to do this, we use of the coupling results of \cite{Dedecker2004a,Dedecker2005} under $\tau$-dependence assumptions to use \citet{VdV}'s results of tightness and asymptotic equicontinuity (under independence) together with the fidis convergence of the EPCF. 

This paper is organized as follows. In Section 2, we recall basic definitions and notations: cluster functionals, the triangular arrays (or normalized random variables excesses) and examples. Then, we give the definition of cluster functionals empirical processes and close the section with a simple version of the CLT of those empirical processes. In Section 3 we begin with the fidis convergence of the EPCF, followed by the conditions to obtain asymptotic tightness and asymptotic equicontinuity of thoses processes to obtain uniform convergence. We close this section with an application: block estimator of the extremal index. In Section 4 we develop a example similar to (\ref{5}) for the multidimensional case for the case of the AR(1)-inputs (\ref{AR}). Also a simulation study for the extremal index to demonstrate the accuracy of this technique. The $\tau$-weak dependence with some examples and the proofs are reported in Appendix. 
\section{Basic definitions and notations}
To define the empirical processes of cluster functionals it is necessary to consider two important ingredients: the cluster functionals and the excesses over high thresholds.

 \subsection{Cluster Functionals}
Let $(E,\mc{E})$ be a measurable subspace of $(\RR^d,\mc{B}(\RR^d))$ for some $d\³1$ such that $0\in E$. Following the deterministic definition of  \citet{Drees2010}\footnote{This definition is given by  \cite{Yun2000} and \cite{Segers2003}, for the real case}, we consider the set of $E$-valued sequences of finite length, {\it i.e.}, 
$$E_\cup:=\{(x_1, \ldots, x_r): x_i\in E  \quad\forall i=1,\ldots, r; \quad\forall r\in\NN \},$$
equipped with the $\sigma$-field $\mc{E}_\cup$ induced by Borel-$\sigma$-fields  on $E^r$, for $r\in\NN$. Let $x\in E_\cup$, then we can write $x=(x_1, \ldots, x_r)$ for some $r\in\NN$. The \textbf{core}\footnote{ Note that the core also considers the null values that exist between the non-null values. For example. $(0,1,2,4,0,3,0,1,0,0)^c=(1,2,4,0,3,0,1)$, which is the smaller sub-block of $x=(0,1,2,4,0,3,0,1,0,0)$ which contains all non-null values as well as the null values between them.} $x^c\in E_\cup$ of  $x$ is defined by 
 \[x^c:=\left\{\begin{array}{cc}
 (x_{r_I}, x_{r_I+1}, \ldots, x_{r_S}), & \mbox{ if  $x\­ 0_r$ (the null element in $E^r$)}    \\ \\
 0,  & \mbox{  otherwise}
\end{array}\right.\]
where $r_I: =\min\{i\in \{1,\ldots,r\}: x_i \­ 0\}$ (first non-null value of the block $x$) and $r_{S}: =\max\{i\in \{1,\ldots,r\}: x_i \­ 0\}$ (last non-null value of the block $x$). A \textbf{cluster functional} is a measurable map $f: (E_\cup, \mc{E}_\cup)\lfled (\RR, \mc{B}(\RR))$ such that 
 \begin{align}\label{fc}
 f(x)=f(x^c),\qquad\mbox{ for all }x\in E_\cup,\ \ \mbox{ and } \ f(0_r)=0\  (\forall r\³1).
 \end{align}
  
Under the properties (\ref{fc}), it is easy to build a large amount of examples of cluster functionals. Nevertheless, the typical examples used to build estimators through these cluster functionals are functionals of the type: 
\begin{align}\label{combi}
 f(x_1, \ldots, x_r)=\sum_{i=1}^r \phi(x_i),
\end{align}
where $\phi: E\lfled \RR$ is a measurable function such that $\phi(0)=0$. Generally speaking, these functions $\phi$  are  indicator functions (or functions which are product of another measurable function $H:E \lfled \RR$  with an indicator function). Another classic example is the component-wise maximum of a cluster:
 \begin{align}\label{combi2}
 f(x_1, \ldots, x_r)=\max_{1\²i\²r} x_i, 
\end{align}
for $E=[0,\infty)$.

Under the set $E=[0,\infty)$, two particular examples that we shall only mention here with a view to motivating further work, are the following functionals:
\begin{itemize}
\item Balanced periods  at $u>0$, 
 
\[f(x_1,\ldots,x_r)=\left\{\begin{array}{ll}
\mb{1}{\left\{\sum_{i=1}^r (x_i - u)\mb{1}\{x_i>0\}=0\right\}}, & \mbox{if  $x_i>0$ for some $i=1,\ldots,r$}     \\
 0,  & \mbox{if $x_i=0$ for all $i=1,\ldots, r.$}
\end{array}\right.\]

\item Maximum sum (greater than the level $u>0$) of  consecutive excesses,
$$f(x_1,\ldots,x_r)=\max_{1\²p<q\²r}\left(\sum_{y \in H_{p,q}}y\right) \mb{1}\left\{\sum_{y\in H_{p,q}} y>u\right\},$$
with the notation: $H_{p,q}=\{x_p,x_{p+1},\ldots,x_q\}\subseteq\{x_1,\ldots,x_r\}$ such that $y>0$, $\forall y\in H_{p,q}$, where  $p,q\in\{1,2,\ldots,r\}$.
\end{itemize}

\subsection{The excesses and their normalizations}
Without going yet into formalities, first let us consider the following examples that motivate the use of the triangular arrays $(X_{n,i})_{1\²i\²n, n\in\NN}$ throughout this work. 

1.-  Let $(X_i)_{i\in\NN}$ be a real-valued stationary random process with marginal cumulative distribution function $F$ and let $(u_n)_{n\in\NN}$ be a non-decreasing sequence of thresholds such that $u_n \uparrow x_F$, where $$x_F=\sup\{x\in\RR: F(x)<1\},\qquad v_n=\PP\{X_1>u_n\} \underset{n\to\infty}{\lfled} 0.$$ 
If we want to study the process $X'_{n,i}= X_i-u_n | X_1>u_n$, first observe that the tail distribution function of $X_i$ may be  asymptotically degenerated as $n\to\infty$, which means that there exists a point $a\in\RR$ such that $$\overline{P}_n(x)=\PP\{X_1-u_n>x|X_1>u_n\}\underset{n\to\infty}{\lfled} \1\{x\²a\}.$$  
However, if $F$ belongs to the  domain of attraction of some extreme-value distribution, then by a result in  \citep{Pickands1975}, there exists $\gamma\in\RR$ and a sequence of positive constants $(a_n)_{n\in\NN}$ (depending on the sequence $u_n$) such that
\[P_n(x)=\PP\{X_{n,1}>x|X_1>u_n\} \underset{n\to\infty}{\lfled} \left\{\begin{array}{ll}
(1+\gamma x)^{-1/\gamma}_+, & \mbox{if  $\gamma\­0$}     \\
\e^{-x},  & \mbox{if $\gamma=0$}
\end{array}\right.\]
locally uniform in $(0,\infty)$, where 
\begin{align}\label{N1}
X_{n,i}=\left(\frac{X_i-u_n}{a_n}\right)_+:= \max\left\{ \frac{X_i-u_n}{a_n} , 0 \right\}, \qquad \text{ for }1\leqslant i \leqslant n;
\end{align}
are the normalized excesses of $X_i$ over $u_n$.

2.- For the multidimensional case we may consider the following example. For $d\³1$, let  $(\mathbf{X}_i)_{i\in\NN}$ be a $\RR^d$-valued random process such that $\mathbf{X}_i=(X^{(1)}_{i},X^{(2)}_{i},\ldots ,X^{(d)}_{i})$ admits coordinates with the same marginal distribution, then in this case, a standardization of $\mathbf{X}_i$ is:
\begin{align}\label{3}
X_{n,i}=\left(\left(\frac{X^{(1)}_i-u_n}{a_n}\right)_+, \ \left(\frac{X^{(2)}_i-u_n}{a_n}\right)_+,\ \ldots, \ \left(\frac{X^{(d)}_i-u_n}{a_n}\right)_+\right),
\end{align}
where $(u_n)_{n\in\NN}$ and $(a_n)_{n\in\NN}$ are defined as in eqn. (\ref{N1}). Here, $X_{n,i}$ is  the vector of normalized excesses over the threshold $u_n$  for each coordinate.

Another interesting example is the normalization of $d$ consecutive excesses of real-valued random variables $(X_i)_{i\in\NN}$,  {\it i.e.}
 \begin{align}\label{4}
X_{n,i}=\left(\left(\frac{X_i-u_n}{a_n}\right)_+, \left(\frac{X_{i+1}-u_n}{a_n}\right)_+,\ldots, \left(\frac{X_{i+d-1}-u_n}{a_n}\right)_+\right).
\end{align}
This example is given in Section 3 - \citep{Drees2010}. Observe that this example is a particular case of the example of eqn.  (\ref{3}) which corresponds to $X^{(j)}_i=X_{i+j-1}$ for $1\le j\le d$.)\\ 
Notice that this example brings information on the extremal dependence structure. Some applications of this standardization could be:\\
{\it  (i) $d$ consecutive days of rain are observed  in a given city, such that the volume of precipitations may be larger  than the volume of water that can be drained (through sewers, soil, rivers, etc.),
 \\
  (ii) $d$ very large claims are reported to an insurance company in a very small time interval  (with respect to typical cases) which this can be a  risk  with respect to the response capacity of the insurance company, and 
  \\
  (iii) $d$ consecutive days of low temperatures observed in a given city, such that the power consumption (due to heating, etc.)  endangers the response capacity of the company in charge of the energy distribution.}

In a general way, let $(E,\mc{E})$ be a measurable subspace of $(\RR^d,\mc{B}(\RR^d))$ for some $d\³1$ such that $0\in E$. We denote by $(X_{n,i})_{1\²i\²n}$ as the $E$-valued row-wise stationary standardized random variables, defined on some probability space $(\Omega, \mc{A}, \PP)$, which are built from a stationary random process $(X_i)_{i\in\ZZ}$, in a way such that the standardization $X_{n,i}$ maps all "non-extreme" values to zero. Additionally, it should satisfy that the sequence of conditional distributions of $X_{n,1}$ given that $X_{n,1}$ belongs of the failure set $S \subseteq E\setminus \{0\}$ ({\it i.e. }$P_n(\cdot| S):=\PP\{X_{n,1}\in \cdot  |X_{n,1}\in S\}$),  converge weakly to some non-degenerate limit.

\subsection{Empirical Processes of Cluster Functionals}
Now, we want to apply cluster functionals to blocks of $E$-valued random variables excesses over a determined thresholds sequence and to define the empirical process indexed by these functionals. 

In order to do that, first let us consider a row-wise stationary $E$-valued triangular array $(X_{n,i})_{1\²i\²n, n\in \NN}$, defined on some probability space $(\Omega, \mc{A}, \PP)$.

Let $Y_{n,j}$ be the $j$-th block of $r_n$ consecutive values of the $n$-th row of $(X_{n,i})$. That is, we have $m_n:=[n/r_n]=\max\{j\in \NN: j \² n/r_n\} $ blocks 
\begin{align}\label{bloc}
Y_{n,j}:=(X_{n,i})_{(j-1)r_n+1\²i\² jr_n}
\end{align} of length $r_n$,   with $1\² j \² m_n$.  In order to simplify future notations, since $(X_{n,i})_{1\²i\²n}$ is stationary for each $n$, then we can denote by $Y_n$  to the "generic block" such that $Y_n \overset{\mc{D}}{=}Y_{n,1}$.\\
Let $\mc{F}$ be a class of  cluster functionals. The \textbf{"empirical process $Z_n$ of cluster functionals''} in $\mc{F}$, is the process $(Z_n(f))_{f\in\mc{F}}$ defined by
\begin{align}\label{PE}
Z_n(f):= \frac{1}{\sqrt{n v_n}}\sum_{j=1}^{m_n}(f(Y_{n,j})-\EE f(Y_{n,j})),
\end{align}
where  $v_n:=\PP\{X_{n,1} \in S\}$ and $S$ is the failure set. 

In order to begin approaching the convergence in fidis of the EPCF (\ref{PE}), observe that if the blocks $(Y_{n,j})_{1\²j\²m_n, n\in\NN}$ are independents and if we take in account the following essential convergence assumptions:
 \begin{itemize}
\item[\textbf{(C.1)}] $\EE\left[\left(f(Y_{n})-\EE f(Y_{n})\right)^2\mb{1}\left\{|f(Y_{n})-\EE f(Y_{n})|>\epsilon\sqrt{n v_n}\right\}\right]=o(r_n v_n)$, \\for all $\epsilon>0$, and  for all $f\in\mc{F}$.
\item[\textbf{(C.2)}] $(r_n v_n)^{-1}\Cov\left(f(Y_{n}),g(Y_{n})\right)\lfled c(f,g)$, for all $f,g\in \mc{F}$,
\end{itemize}
with  $r_n \ll v_n^{-1} \ll n$, then the fidis of the empirical process  $(Z_n(f))_{f\in\mc{F}}$ of cluster functionals converge to the fidis of a Gaussian process $(Z(f))_{f\in\mc{F}}$ with the covariance function $c$.  

\citet{Drees2010} have proved CLTs for the process (\ref{PE}). In particular, they have proved CLTs in fidis for this process by using the  Bernstein blocks technique together with a  $\beta$-mixing coupling condition to boil down convergence to convergence of sums over i.i.d. blocks through \citet{Eberlein}'s technique involving the metric of total variation.  Moreover, \citet{Drees2010} extend the results to the uniform convergence by adding \citet{VdV}'s tightness criteria and asymptotic equicontinuity conditions to the results in fidis that they had obtained. 

We aim at  extending  their CLT's for  the empirical process $(Z_n(f))_{f\in\mc{F}}$, since the family of mixing processes is still very restrictive. One particular and really simple example of a non-mixing process is the AR(1)-process (\ref{AR}). We derive some results as in \citep{Drees2010} and some applications as in \citep{Drees2011} under much  weaker dependence conditions including {\it eg.} this example. 
\\
The $\tau-$weak dependence introduced by \citet{Dedecker2004a} holds  for the Example in eqn.  (\ref{AR}), as well as  more generally for  Bernoulli shifts processes and Markov chains.
\section{Central limit theorems for cluster functionals}

\subsection{Fidis convergence}
First we give the convergence in fidis. In this case, the technique used to prove the convergence of the empirical process  (\ref{PE}) is also the Bernstein blocks technique. In order to do this, we need to extract from each block $Y_{n,j}$ of length $r_n$ a sub-block of  length $l_n$, in such a way that $l_n=o(r_n)$. Then we use the remaining sub-blocks, separated by $l_n$ variables, combined with  convenient conditions of $\tau$-dependence to couple independent blocks to the original blocks (originally dependent), and thus obtain the CLTs through classic tools. 

For this, it is necessary for the triangular array $(X_{n,i})_{1\²i\²n, n\in\NN}$ to satisfy the relation:
\begin{itemize}
\item[\textbf{(D.1)}]  $\tau_{1,n}(l_n)=o(r_n^{-1})$,
\end{itemize}
such that
\begin{itemize}
\item[\textbf{(B.1)}]  $l_n \ll r_n \ll v_n^{-1} \ll n$ where $l_n\lfled \infty$ and $n v_n\lfled \infty$ as $n\to\infty$, and
\item[\textbf{(B.2)}] $\EE\left(\|X_{n,1}\| | X_{n,1}\­0\right)<\infty.$
\end{itemize}

On the other hand, we must not forget the influence of the small blocks extracted with length $l_n$ over the empirical process  (\ref{PE}). In this case, in order to introduce assumptions over these small blocks, it is necessary to consider the following notations, which we will also use throughout the rest of this paper. 

\begin{nota}\rm Let $x=(x_1, x_2, \ldots, x_r)$. We will use the notation $x^{[l:k]}$ as follows
\[x^{[l:k]}=\left\{\begin{array}{ll}
0, & \mbox{if $r<l$} ,\\ 
 (x_l,\ldots, x_k), & \mbox{if  $1\²l\²k \² r$},   \\ 
 Y, & \mbox{if $k > r$.}
\end{array}\right.\]
and $x^{[k]}:=x^{[1:k]}$. Moreover, if $f\in\mc{F}$ is a cluster functional, then we denote
\begin{align}\label{Delta}
\Delta_n(f):=f(Y_n)-f(Y_n^{[r_n-l_n]}),
\end{align}
where $r_n$ is the length of the block $Y_n$ and $l_n$ is such that $l_n=o(r_n)$.
\end{nota}

The following assumption guarantees that the extraction of the small blocks does not disturb the result of the convergence in fidis (if that is the case) of the EPCF $(Z_n(f))_{f\in\mc{F}}$. 
\begin{itemize}
\item[\textbf{(C.3)}] For all $f\in\mc{F}$, 
\begin{align} \EE|\Delta_n(f)-\EE \Delta_n(f)|^2\mb{1}\left\{|\Delta_n(f)-\EE \Delta_n(f)|\²\sqrt{n v_n}\right\}=&o(r_n v_n)
\nonumber\\ \pp\left\{|\Delta_n(f)-\EE \Delta_n(f)|>\sqrt{n v_n}\right\}=&o(r_n/v_n).\nonumber
\end{align}
\end{itemize}
\begin{theo} Suppose that (B.1)-(B.2),  (C.1)-(C.3) and (D.1) hold. Then the fidis of the cluster functionals empirical process $(Z_n(f))_{f\in\mc{F}}$ converge to the fidis of a centered Gaussian process $(Z(f))_{f\in\mc{F}}$ with covariance function $c$ defined in assumption (C.2).
\end{theo}
\begin{rem}\rm Note that Assumptions (C.1)  and (C.3) are difficult to check in general, for that reason, consider the following (more restrictive but easier to verify) alternatives conditions:
\begin{itemize}
\item[\textbf{(A.1)}] $\Var(\Delta_n(f))=o(r_n v_n)$ for all $f\in\mc{F}$.
\item[\textbf{(A.3)}] $\EE (f(Y_n))^{2+\delta}= O(r_n v_n)$ for some $\delta>0$ and for all $f \in\mc{F}$. 
\end{itemize}
\end{rem}
\begin{lem} The conditions (A.1) and (A.3) implies the conditions (C.1) and (C.3), respectively.
\end{lem}
For the proof of this lemma see Lemma 5.2 in \citep{Drees2010}.

Generally, the (C.2) convergence can be easily verified. However, in some situations it could be difficult to give the limit $c$ in an explicit mode. In this sense, it is common to use the "tail chain" associated to the original process $(X_i)_{i\in\NN}$. This terminology was used by \citet{Perfekt1994}, but generalized by \citet{Yun1998} to make explicit representations of the extremal index of a higher-order stationary Markov chain. \citet{Segers2003} generalizes this result to stationary sequences with suitable conditions. 

We will follow \citet{Segers2003} rationale to give an explicit representation of the covariance function $c$ through the tail chains. The following proposition (which is a similar result to \citet{Segers2003}'s Theorems 1 and 3) provide some conditions, sufficient to verify (C.2) and which in some situations are easier to prove. The alternative expression to the covariance function c (defined in (C.2)), is shown below in Corollary 1.

In order to carry this out, it is necessary to consider the following assumption: 
\begin{itemize}
\item[\textbf{(C.2')}] There is a sequence $W=(W_i)_{i\³1}$  of $E$-valued random variables such that, for all $k\in\NN$, the joint conditional distribution $$P_{(X_{n,i}, \mb{1}\{X_{n,i}=0\})_{1\²i\²k}| X_{n,1}\­0}$$ converges weakly to $P_{(W_i, \mb{1}\{W_i=0\})}$, and for all $f\in\mc{F}$ are a.s. continuous with respect to the distribution of $W^{[k]}=(W_1, \ldots, W_k)$ and $W^{[2:k]}=(W_2, \ldots, W_k)$ for all $k$, that is, 
$$\PP\{W^{[2:k]}\in D_{f, k-1}, W_i=0,\ \forall i>k \}=\PP\{W^{[k]}\in D_{f, k}, W_i=0,\ \forall i>k \}=0$$
where we denote by $D_{f,k}$  the set of discontinuities of $f|_{E^k}$.
\end{itemize}
\begin{rem}\rm
The existence of such sequence  $W$ is guaranteed in particular from Theorem 2 in \citep{Segers2003} with  $E=\RR$ and the normalization (\ref{N1}). There, Segers has shown that if $$\PP_{((X_{n,i})_{1\²i\²k}| X_1>u_n)}\underset{n\to\infty}{\lfled} -\log G_k,$$  where $G_k$ is some $k-$dimensional extreme value distribution for all $k\in\NN$, then there exists such "{\it tail chain}" $W=(W_i)_{i\in\NN}$ such that 
\begin{align}\label{tailchains}
P_{((X_{n,i}, \mb{1}\{X_{n,i}=0\})_{1\²i\²k}| X_1>u_n)} \overset{w}{\underset{_{n\to \infty}}{\lfled}} P_{(W_i, \mb{1}\{W_i=0\})_{1\²i\²k}},
\end{align} 
for all $k\in\NN$. 
\end{rem}
\begin{prop}  Suppose that the r.v's $(X_{n,i})_{1\²i\²n, n\in\NN}$ satisfies the following condition:
\begin{itemize}
\item[\textbf{(D.2)}] There exists $p,q\³1$ with  $p^{-1}+q^{-1}=1$ and $\alpha>0$, such that 
\begin{align}\label{Tau}
\lim_{l\to \infty}\limsup_{n\to \infty} \dfrac{r_n \tau_{p,n}(l)}{v_n^{1/p+\alpha}}=0\qquad \text{and   }\qquad r_n v_n^{1/q} \underset{n\to\infty}{\lfled} 0
\end{align}
\end{itemize}
Then, 
$$\EE\left[f(Y_n)| Y_n\­0\right]=\theta_n^{-1} \EE\left[f(Y_{n,1})-f(Y_{n,1}^{[2:r_n]}) | X_{n,1}\­0\right] + o(1),$$
where $o(1)$ converges to $0$ as $n\to \infty$ uniformly for all bounded cluster functionals $f\in\mc{F}$, and 
\begin{align}\label{index}
\theta_n:= \dfrac{\PP\{Y_n\­ 0\}}{r_n v_n}=\PP\{ Y_{n,1}^{[2: r_n]}=0| X_{n,1}\­ 0\}(1+o(1)).
\end{align}
Additionally, if the assumption (C.2') is satisfied, then:
$$
\begin{array}{cll}
m_W &:=&\sup\{i\³1: W_i\­0\}< \infty,
\\ \theta_n& \underset{_{n\to \infty}}{\lfled}&\theta:=  \PP\{ W_i=0, \forall i\³ 2\}=\PP\{ m_W=1\}>0,
\\ P_{f(Y_n)| Y_n\­ 0}& \overset{w}{\underset{_{n\to \infty}}{\lfled}}& \dfrac{1}{\theta} \left(  \PP\{f(W)\in \cdot \}- \PP\{f(W^{[2:\infty]})\in \cdot , m_W\³ 2\}\right).
\end{array}
$$
\end{prop}
\begin{coro} Suppose  that $$\mc{F}=\{f|  (f(Y_n)^2)_{n\in\NN} \text{ is uniformly integrable under } P(\cdot )/r_nv_n\}.$$ Assume that (B.1), (B.2), (C.2'), (C.3) and (D.2) hold. Then the fidis of the cluster functionals empirical process $(Z_n(f))_{f\in\mc{F}}$ converge to the fidis of a centered Gaussian process $(Z(f))_{f\in\mc{F}}$ with covariance function $c$ defined by
\begin{align}\label{cc}
c(f,g)=\EE\left[(fg)(W)- (fg)(W^{[2:\infty]})\right].
\end{align}
\end{coro}
There are many cases in which $\|f\|_\infty=\sup_{x\in E_{\cup}}|f(x)|<\infty$, for all $f\in\mc{F}$. Under this condition, it is clear that the conditions (C.1) and (C.3) are satisfied. Therefore, it is important to note the following corollary. 
\begin{coro} Suppose that (B.1), (B.2), (C.2') and (D.2) are satisfied. Then, if $\|f\|_\infty=\sup_{x\in E_{\cup}}|f(x)|<\infty$ for all $f\in\mc{F}$, the fidis of the cluster functionals empirical process $(Z_n(f))_{f\in\mc{F}}$ converges to the fidis of a centered Gaussian process $(Z(f))_{f\in\mc{F}}$ with covariance function $c$ defined by (\ref{cc}).
\end{coro}
\subsection{Uniform convergence} To prove uniform convergence, we use either asymptotic tightness of $Z_n$ in the space $\ell^\infty(\mc{F})$, or asymptotic equicontinuity conditions, by some results in \S\, 2.11 of \citet{VdV}. Those results need independence therefore a argument of coupling for the blocks $(Y_{n,j})_{1\²j\²m_n, n\in\NN}$ is also used here. 

\subsubsection{Asymptotic tightness}
\begin{defi}  The sequence $(Z_n)_{n\in\NN}$ is asymptotically tight if for every $\epsilon>0$ there exists a compact set $K\in\ell ^\infty(\mc{F})$ such that
$$\limsup_{n\lfled \infty} \PP^*(Z_n\notin K^\delta)< \epsilon, \quad \text{ for every }\delta>0,$$
where $K^\delta=\{f\in\ell^\infty(\mc{F}): d_{\mc{F}}(f,K)<\delta\}$ is the "$\delta-$enlargement" around $K$ and $\PP^*$ denotes the outer probability.
\end{defi}
\begin{defi}  The bracketing number $N_{[\cdot]}(\epsilon, \mc{F},L_2^n)$ is defined as the smallest number $N_\epsilon$ such that for each $n\in\NN$ there exits a partition $(\mc{F}_{n,k}^\epsilon)_{1\²k\²N_\epsilon}$ of $\mc{F}$ such that 
$$\EE^*\sup_{f,g\in \mc{F}_{n,k}^\epsilon}\left(f(Y_n)-g(Y_n)\right)^2 \² \epsilon^2 r_n v_n, \hspace{0,5cm} \mbox{ for } \quad 1\² k\² N_\epsilon,$$
where $\EE^*$ denotes the outer expectation. 
\end{defi}
 In order to use  Theorem 2.11.9 in \citep{VdV} we need:
\begin{itemize}
\item[\textbf{(T.1)}]  The set $\mc{F}$ of cluster functionals is such that for each $f\in\mc{F}$ the expression $\EE f^2(Y_n)$ is finite for all $n\in\NN$ and such that the envelope function  satisfies:
$$F(x):=\sup_{f\in\mc{F}} |f(x)|< \infty, \qquad \forall x\in E_\cup.$$
\item[\textbf{(T.2)}]  $\EE^*\left(F(Y_n)\mb{1}\{F(Y_n)>\epsilon \sqrt{n v_n}\}\right)=o(r_n \sqrt{v_n/n})$ for all $\epsilon > 0$.
\end{itemize}
Note that for a sequence of monotonically increasing positive functions $(h_n(\delta))_{n\³1}$ the convergence of $h_n(\delta_n)$ to zero $\forall \delta_n \downarrow 0$ is equivalent to $$\lim_{\delta\downarrow 0}\limsup_{n\lfled \infty} h_n(\delta)=0,$$  thus the Assumptions 2 and 3 of Theorem 2.11.9 from \citep{VdV} are  reformulated as follows: 
\begin{itemize}
\item[\textbf{(T.3)}]  There exists a semi-metric $\rho$ on $\mc{F}$ such that $\mc{F}$ is totally bounded with respect to (w.r.t.) $\rho$ and 
$$\lim_{\delta \downarrow 0} \limsup_{n\lfled \infty} \sup_{f,g\in\mc{F}: \rho(f,g)<\delta} \dfrac{1}{r_n v_n}\EE \left(f(Y_n)-g(Y_n)\right)^2=0.$$ 
\item[\textbf{(T.4)}]  $$\lim_{\delta \downarrow 0} \limsup_{n\lfled\infty} \int_0^\delta \sqrt{\log N_{[\cdot]}(\epsilon, \mc{F},L_2^n)}d\epsilon=0.$$
\end{itemize} 
\begin{theo} Suppose that (B.1), (B.2), (D.1) hold and that (T.1) - (T.4)  are satisfied. Then the empirical process $(Z_n)_{n\in\NN}$ is asymptotically tight in $\ell^\infty(\mc{F})$. Moreover, if  the assumptions (C.1)-(C.3) hold, then  $Z_n$ converges to a centered Gaussian process $Z$ with covariance function $c$ in (C.2).
\end{theo}
\subsubsection{Asymptotic equicontinuity}
\begin{defi} The sequence $(Z_n)_{n\in\NN}$ is asymptotically equicontinuous w.r.t. a semi-metric $\rho$ if for any $\epsilon>0$ and $\eta>0$  there exists some $\delta>0$ such that: 
$$\limsup_{n\lfled \infty} \PP^*\left(\sup_{f,g\in\mc{F}:\; \rho(f,g)<\delta}|Z_n(f)-Z_n(g)|>\epsilon\right)< \eta.$$\end{defi}
We use Theorem 2.11.1 in \citep{VdV} to prove asymptotic equicontinuity. In order to do this, we need to define a semi-metric $\rho_n$ on $\mc{F}$ as follows. Let $(Y^*_{n,j})_{1\²j\²m_n}$ be the independent copies of  $(Y_{n,j})_{1\²j\²m_n}$.\\ We define $\rho_n$ as:
\begin{align}\label{Semi}
\rho_n(f,g):=\sqrt{\dfrac{1}{n v_n}\sum_{j=1}^{m_n}(f(Y_{n,j}^*)-g(Y_{n,j}^*))^2}.
\end{align}
Moreover,  we denote by $N(\epsilon, \mc{F},\rho)$ the "covering number",  the minimum number of balls (w.r.t. the semi-metric $\rho$)  with radius $\epsilon>0$ necessary to cover $\mc{F}$. In this way, we can add to the list of assumptions the following two:
\begin{itemize}
\item[\textbf{(T.4')}]  For $k=1,2$ the map
$$(Y^*_{n,1},\ldots, Y^*_{n,[m_n/2]}) \longmapsto \sup_{f,g\in\mc{F}: \rho(f,g)<\delta}\sum_{j=1}^{[m_n/2]} e_j \left(f(Y_{n,j}^*)-g(Y_{n,j}^*)\right)^k$$
is measurable for each $\delta>0$, each vector $(e_1,\ldots,e_{[m_n/2]})\in\{-1,0,1\}^{[m_n/2]}$ and each $n\in\NN$.
\item[\textbf{(T.5)}]  $$\lim_{\delta \downarrow 0} \limsup_{n\lfled\infty} \PP^*\left(\int_0^\delta \sqrt{\log N(\epsilon, \mc{F},\rho_n)}d\epsilon>\eta\right)=0\text{,  } \hspace{0.5cm}\forall \eta >0.$$
\end{itemize}
\begin{theo} Suppose that (B.1), (B.2), (D.1) hold and that  (T.1)-(T.3), (T.4') and (T.5) are satisfied. Then the empirical process $(Z_n)_{n\in\NN}$ is asymptotically equicontinuous. Moreover if  the assumptions (C.1)-(C.3) hold, then  $Z_n$ converges to a centered Gaussian process $Z$ with covariance function $c$ in (C.2).
\end{theo}
\subsection{Application: blocks estimator of the extremal index}
Let  $(X_i)_{i\in\NN}$ be a real stationary time series with distribution function $F$. Now consider the index defined in (\ref{index}) with the extreme normalization (\ref{N1}) and $u_n:=F^{\leftarrow}(1-v_n t)$, for all $t\in [0,1]$, {\it i.e.}
\begin{align}\label{index2}
\theta_{n,t}:= \dfrac{\PP \{ Y_n \­ 0\}}{r_n v_n t}= \dfrac{\PP \{ \max_{1\² i\² r_n} X_i > u_n\} }{r_n v_n t},  \quad \text{ with } t\in (0, 1].
\end{align}
Note that if $r_n$ satisfies (B.1) condition and if (D.2) holds, then by using Proposition 1, there exists a number ({\it extremal index}) $\theta\in (0,1]$  such that 
\begin{align}\label{convergen}
\theta_{n,t}\underset{_{n\to \infty}}{\lfled} \theta  \quad \text{ uniformly for } t\in (0, 1].
\end{align}
Given the convergence (\ref{convergen}), Drees has suggested in his paper \cite{Drees2011} to estimate $\theta$ replacing the unknown probability $\PP\{ \max_{1\²i\²r_n} X_i>u_n\}$ and the unknown expectation $r_n v_n t=\EE\left[\sum_{i=1}^{r_n} \mb{1}\{ X_i> u_n\} \right]$ by a empirical expression for $\theta_{n,t}$:
\begin{align}\label{t.empi}
\h{\theta}_{n,t}:=\dfrac{\sum_{j=1}^{m_n}\mb{1}\{\max_{(j-1)r_n< i \² jr_n}X_i >u_n\}}{\sum_{j=1}^{m_n}\sum_{i=(j-1)r_n+1}^{jr_n}\mb{1}\{ X_i> u_n\}},
\end{align}
where $m_n=[n/r_n]$ such that $1 \ll r_n \ll v_n^{-1} \ll n$ but $n v_n \lfled \infty$. Thus, such estimator (\ref{t.empi}) (called {\it blocks estimator of the extremal index}) can be expressed in terms of two empirical processes of cluster functionals $(Z_n(f_t), Z_n(g_t))_{0\²t\² 1}$ defined in (\ref{PE}). For this, suppose without loss of generality that the random variables $(X_i)_{1\²i\²n}$ are uniformly distributed on $[0,1]$ (otherwise, just consider the transformation $U_i=F(X_i)$, $1\²i\²n$, where $F$ is the distribution function of $X_1$, see \cite{Drees2011}). Then, with the normalization (\ref{N1}) such that $a_n=v_n=1-u_n$ and the blocks $(Y_{n,j})_{1\²j\²m_n}$ defined in (\ref{bloc}), we have that
\begin{align}\label{est.block}
\h{\theta}_{n,t}=\dfrac{m_n^{-1}\sum_{j=1}^{m_n} f_t(Y_{n,j})}{m_n^{-1}\sum_{j=1}^{m_n} g_t(Y_{n,j})}= \dfrac{\EE f_t(Y_{n,1}) +(n v_n)^{1/2} m_n^{-1} Z_n(f_t)}{\EE g_t(Y_{n,1})+(n v_n)^{1/2} m_n^{-1} Z_n(g_t)},
\end{align}
where
\begin{align}\label{F1}
f_t(x_1, \ldots, x_r)&:= \mb{1}\{ \max_{1\²i\² r} x_i>1-t\}
\\ g_t(x_1, \ldots, x_r)&:= \sum_{i=1}^r \mb{1}\{ x_i>1-t\}.
\end{align}
 For this particular case, we consider the following assumptions: 
 \begin{itemize}
\item[\textbf{(C.2.1)}] $(r_n v_n)^{-1} \Cov (g_s(Y_n), g_t(Y_n)) \lfled c_g(s,t)$, for all $0\² s,t\²1$.
\item[\textbf{(C.2.2)}] $(r_n v_n)^{-1} \Cov (f_s(Y_n), g_t(Y_n)) \lfled c_{fg}(s,t)$, for all $0\²s,t\²1$.
\item[\textbf{(T)}] For some bounded function $h: (0,1]\lfled \RR$ such that $\lim_{t\to 0}h(t)=0$
$$(r_n v_n)^{-1}\EE\left(f_s(Y_{n,1})-f_t(Y_{n,1}) \right)^2\² h(t-s), \quad \forall 0\²s<t\²1,$$ for all $n$ sufficiently large. 
\end{itemize}
The following are a slight variation of the first two results of \cite{Drees2011}, in the sense that we replace the $\beta$-mixing condition with $\tau$-dependence condition. 
\begin{prop}\
\begin{itemize}
\item [(3.1)] Suppose that (B.1), (B.2) and (D.1) are satisfied. Then $(Z_n(f_t))_{0\²t\²1}$ converges weakly to $Z_f:=(\sqrt{\theta} B_t)_{0\²t\²1}$, where $B$ denote a standard Brownian motion.
\item [(3.2)] If additionally (C.2.1) and (T) are satisfied and $r_n=o(\sqrt{nv_n})$, then the sequence of processes $(Z_n(g_t))_{0\²t\² 1}$ converges weakly to a centered Gaussian process $(Z_g(t))_{0\²t\²1}$ with covariance function $c_g$.
\item [(3.3)] Under all the hypothesis of (3.1) and (3.2), if moreover (C.2.2) holds, then $(Z_n(f_t),Z_n(g_t))_{0\²t\²1}$ converge weakly to $(Z_f(t),Z_g(t))_{0\²t\²1}$ with 
\begin{align}
\Cov(Z_f(s),Z_f(t))&=\theta( s \land t),
\nonumber\\ \Cov(Z_g(s),Z_g(t))&= c_g(s,t),
\nonumber\\ \Cov(Z_f(s),Z_g(t))&= c_{fg}(s,t),  \quad 0\²s,t\²1.
\end{align}
\end{itemize}
\end{prop}

Using the same argument in Remark 2, we can find explicit expressions for the covariance functions $c_g$ and $c_{fg}$ as functions of the "tail chains" of $(X_i)_{i\in\NN}$. This is,  if for every $k \in \NN$ the distribution function of $(X_1, \ldots, ,X_k)$ belongs to the domain of attraction of an extreme-value distribution, then there exist a sequence $W=(W_i)_{i\in\NN}$ such that (\ref{tailchains}) hold. In such case:
\begin{align}
c_g(s,t)&= s \land t + \sum_{i=1}^{\infty} \left( \PP \{ W_1>1-s, W_{i+1}>1-t\}+ \PP \{ W_1>1-t, W_{i+1}>1-s\} \right)
\nonumber\\ c_{fg}(s,t)&=\left\{\begin{array}{ll} \PP\{W_1>1-t, \max_{j\³ 1} W_j> 1-s \} &  \\ + \sum_{i=1}^\infty \PP\{ W_1>1-s, W_{i+1}>1-t, \max_{j\³ 2} W_j\² 1-s\}, & s<t, \\ t & s\³ t.
\end{array}\right.\nonumber
\end{align}

\begin{coro} Under Proposition 2 - (3.3)'s assumptions, 
\begin{align}\label{gausstheta}
(\sqrt{n v_n} t (\h{\theta}_{n,t}- \theta_{n,t}))_{0< t\² 1}  \overset{w}{\underset{_{n\to \infty}}{\lfled}} Z:= Z_f - \theta Z_g,
\end{align}
where $Z$ is a Gaussian process such that $\EE Z(t)=0$ and 
\begin{align}
\Cov(Z(s), Z(t))=\theta (s \land t - c_{fg}(s,t)- c_{fg}(t,s) + \theta^2 c_g(s,t)).
\end{align}
\end{coro}
\section{Examples and Simulations}
\subsection{AR(1)-process with the functional "number of excesses over $x$"}
We consider the AR(1)-process (\ref{AR})
where $b\³2$ is an integer, $(\xi_i)_{i\in\NN}$ are i.i.d. and uniformly distributed on the set $U(b):=\{0,1,\ldots, b-1\}$. 
\\
It is clear that $X_0$ is uniformly distributed on $[0,1]$. Moreover we define the normalized random variables $(X_{n,i})_{1\²i\²n, n\in\NN}$ as in eqn.  (\ref{4}) with $a_n=v_n=1-u_n$. We set  $(x_1, \ldots, x_d)\²(y_1,\ldots, y_d)$ if and only if $x_i\²y_i$, for all $i=1, \ldots, d$ in case $x,y\in [0,1]^d$. Then
\begin{align}\label{dgn}
\pp&\{X_{n,1}> x| X_{n,1}\­0\}
\nonumber\\ =&\dfrac{1}{b^d \bar{v}_n}\sum_{j_1, \ldots, j_d\in U(b)} \left(\max_{i=1,\ldots, d}\left\{1- b^i +\sum_{s=1}^i b^{s-1}j_s + b^i v_n (1-x_i)\right\}_+ \land 1\right)
\nonumber \\ \underset{_{n\to \infty}}{\lfled} & \max_{i=1,\ldots, d}\{ b^{i-d}(1-x_i)\},
\end{align}
where $\bar{v}_n:=\PP\{ X_{n,1}\­ 0\}\sim \ v_n=\PP\{ X_1> u_n\}\underset{_{n\to \infty}}{\lfled} 0$.

Consider $\mc{F}$  the family of cluster functionals:
\begin{align}\label{func}
\mc{F}=\left\{f_x ,\ {x\in[0,1]^d}\right\}, \quad \mbox{ with }\quad  f_x(x_1, \ldots, x_r)=\sum_{i=1}^r \mb{1}\{x_i>x\}
\end{align}
For this case, we obtain the covariance function $c$ of (C.2):
\begin{align}\label{covar}
c(x,y)=&\min\left( \max_{k=1,\ldots, d}\{ b^k (1-x_k)\}, \max_{k=1, \ldots, d}\{ b^k(1-y_k)\}\right) 
\nonumber\\ &+ \sum_{i=1}^\infty H_{b,i}(x,y) + \sum_{i=1}^\infty H_{b,i}(y,x),
\end{align}
where, for $i\³ d$
\begin{align}\label{cov2}
H_{b,i}(x,y):=\dfrac{1}{b^i}\min\left( \max_{k=1\ldots, d} \{ b^k(1-x_k)\}, \max_{k=1, \ldots, d} \{ b^{k+i} (1-y_k)\}\right) 
\end{align}
and for $1\²i<d$, 
\begin{align}
&H_{b,i}(x,y)
 \nonumber \\ &:=\dfrac{1}{b^i}\min\left( \max_{k=1\ldots, i} \{ b^k(1-x_k)\}, \max_{k=i+1, \ldots, d}\{b^k\min(1-x_k, 1-y_k)\}, \max_{k=d-i, \ldots, d} \{ b^{k+i} (1-y_k)\}\right) \nonumber
\end{align}
Conditions (C.1), (C.3), (T.1)-(T.4) hold for uniformly distributed random variables and for the same family $\mc{F}$, see page 2177 and 2178 in  \cite{Drees2010}. 
Thus,  under  assumption (B.1), setting  $l_n$ and $r_n$ such that $ b^{-l_n}=o( r_n^{-1})$ (see Appendix  A.1 -  Application 1) then the empirical process $(Z_n(x))_{x\in [0,1]^d}$ defined as:
\begin{align}\label{tail}
Z_n(x)&:=\frac{1}{\sqrt{n v_n}}\sum_{i=1}^{r_n m_n} \left( \mb{1}\{X_{n,i}>x\}-\pp\{X_{n,i}>x\}\right)
\nonumber\\ &\sim \frac{1}{\sqrt{n v_n}}\sum_{i=1}^{n} \left( \mb{1}\{X_{n,i}>x\}-\pp\{X_{n,i}>x\}\right)
\end{align}
converges to a centered Gaussian process $Z$ with covariance function (\ref{covar}).
\subsection{Simulation study} The experiment is to estimate the extremal index $\theta$  through the blocks estimator of the extremal index (\ref{est.block}).

Let us consider the AR(1)-process (\ref{AR}). Here, as $X_0$ is uniformly distributed on $[0,1]$ and $X_i= \frac{X_0}{b^i}+\sum_{s=1}^i \frac{\xi_s}{b^{i-s+1}}$ for all $i\³1$, we will take this to obtain a theoretical expression for the index (\ref{index2}) with $u_n=1-v_n t$ for $t\in (0, 1]$:
\begin{align}\label{theoric}
\theta_{n,t}=\dfrac{1}{b^{r_n} r_n v_n t} \sum_{j_1, \ldots, j_{r_n}\in U(b)}\min \left(\max_{i=1,\ldots, r_n}\left\{1-b^i(1-v_n t)+\sum_{s=1}^i b^{s-1} j_s \ \right\}_+,1\right),  
\end{align}
which converges to some $\theta=\theta(b) \in (0, 1)$ if (B.1) is satisfied with $b^{-l_n}=o(r_n^{-1}v_n^\beta)$, for some $\beta>0$.

We will use the advantage of having this theoretical expression (\ref{theoric}) of the extremal index $\h{\theta}_{n,t}$ to compare it to its asymptotic estimations $\h{\theta}_{n,t}$. For this, we simulate $AR(1)-$processes (\ref{AR}) for $b=2,3$ and their blocks estimators (\ref{est.block}) respective with the normalization $(\ref{N1})$ taking $a_n=v_n=1-u_n$  to make the comparison of the results estimated with the theoretic model (\ref{theoric}). 

Let us suppose we have data of a size $M = 600000$ which adjusts appropriately to an $AR(1)$-process (\ref{AR}) for $b=2$ ($b=3$). Because the process is stationary,  we can divide this data into $N = 60$ blocks of length $n = 10^4$. Moreover, we choose a threshold $u_n$ such that $v_n=n^{-1/2}$ and the sub-blocks of lenght $r_n=[\log(n)]$.  In Figure 1 we showed a polygonal curve $(t,\theta_{n=10^4, t})_{t=0.1, 0.2, \ldots, 1}$ (blue curve) of $(t,\theta_{n=10^4, t})_{0\²t\²1}$ and a mean polygonal curve $(t,\h{\theta}_{n=10^4, t})_{t=0.1, 0.2, \ldots, 1}$ (black curve) of $(t,\h{\theta}_{n=10^4, t})_{t=0\²t\²1}$.  Note that for $n=10^4$, the symmetry of the confidence band $(CI_i(t), CI_s(t))_{0\²t\²1}$ with respect to $(\h{\theta}_{n,t})_{0\²t\²1}$ and with a confidence level $1-\alpha=0.95$ (red curves), already shows the gaussian behaviour of the estimator. Furthermore, as expected, the estimated value through the blocks estimator is quite close to the extremal index theoretical (\ref{theoric}), with $n=10^4$. The numerical results are shown in Tables 1 and 2, for the cases $b=2$ and $b=3$, respectively.
\begin{figure}
\centering
\begin{minipage}{.5\columnwidth}
\centering
\includegraphics[width=\columnwidth]{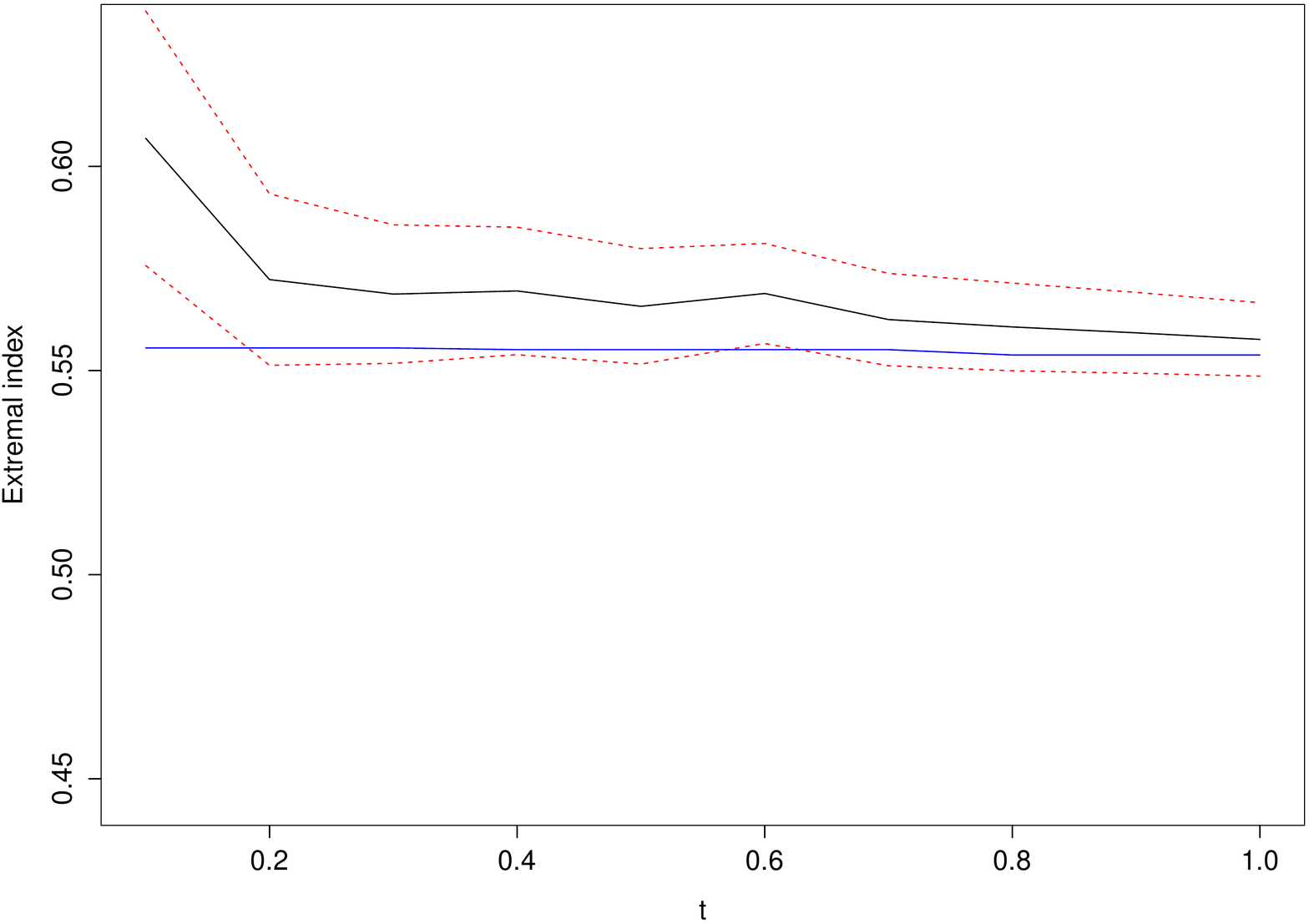}
\label{fig:b2}
\end{minipage}%
\begin{minipage}{.5\columnwidth}
\centering
\includegraphics[width=\columnwidth]{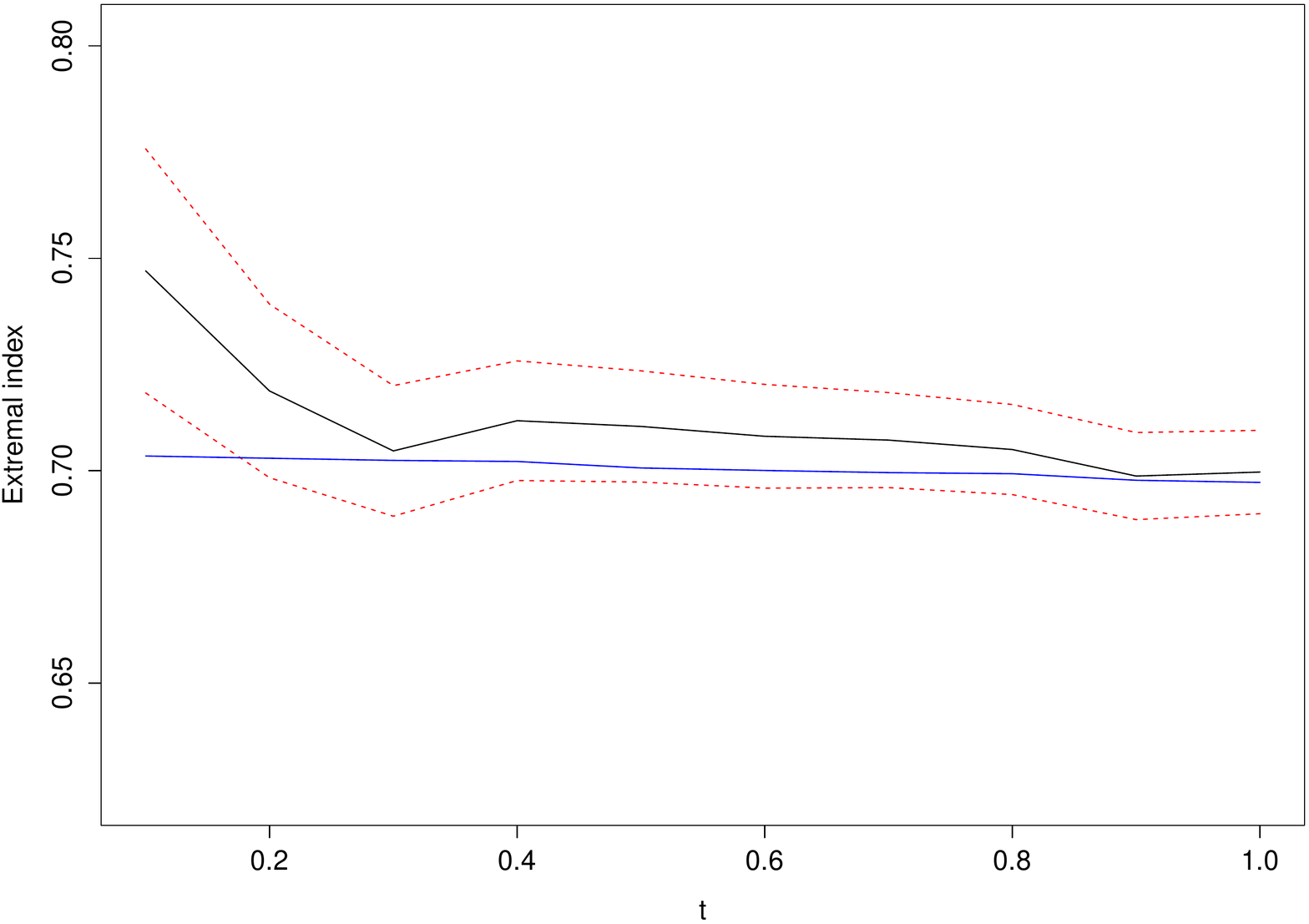}
\label{fig:b3}
\end{minipage}
\caption{\small Extremal index for the AR(1)-input}{\small Left: $\theta_{n=10^4,t}$ is the blue curve, $\h{\theta}_{n=10^4, t}$ is the black curve and the confidence intervals $I_t=(CI_i(t), CI_s(t))$ are the red curves, for the AR(1)-input (\ref{AR}) with $b=2$. Right: the same situation but with $b = 3$}
\end{figure}

\begin{table}[ht]
\centering
\begin{tabular}{rrrrrrrrrrr}
  \hline 
t & 0.1 & 0.2 & 0.3 & 0.4 & 0.5 & 0.6 & 0.7 & 0.8 & 0.9 & 1.0 \\  
  \hline   
   $\theta_{n,t}$  &  0.555 & 0.555 & 0.555 & 0.555 & 0.555 & 0.555 & 0.555 & 0.553 & 0.553 & 0.553 \\ 
  $\h{\theta}_{n,t}$ & 0.606 & 0.572 & 0.568 & 0.569 & 0.565 & 0.568 & 0.562 & 0.560 & 0.559 & 0.557 \\ 
   $CI_s$ & 0.638 & 0.593 & 0.585 & 0.585 & 0.579 & 0.581 & 0.573 & 0.571 & 0.569 & 0.566 \\ 
   $CI_i$ & 0.575 & 0.551 & 0.551 & 0.553 & 0.551 & 0.556 & 0.551 & 0.549 & 0.549 & 0.548 \\   
   \hline
\end{tabular}
\caption{Comparison between the blocks estimation and the theoretical approximation (\ref{theoric}), for the AR(1)-input with $b=2$ and $n=10^4$.}
\end{table}

\begin{table}[ht]
\centering
\begin{tabular}{rrrrrrrrrrr}
  \hline
 t & 0.1 & 0.2 & 0.3 & 0.4 & 0.5 & 0.6 & 0.7 & 0.8 & 0.9 & 1.0 \\  
  \hline
  $\theta_{n,t}$ & 0.703 & 0.702 & 0.702 & 0.702 & 0.700 & 0.700 & 0.699 & 0.699 & 0.697 & 0.697  \\ 
  $\h{\theta}_{n,t}$ & 0.747 & 0.718 & 0.704 & 0.711 & 0.710 & 0.708 & 0.707 & 0.704 & 0.698 & 0.699   \\ 
    $CI_s$ & 0.775 & 0.739 & 0.720 & 0.725 & 0.723 & 0.720 & 0.718 & 0.715 & 0.708 & 0.709  \\ 
  $CI_i$ & 0.718 & 0.698 & 0.689 & 0.697 & 0.697 & 0.695 & 0.696 & 0.694 & 0.688 & 0.689  \\ 
   \hline
\end{tabular}
\caption{Comparison between the blocks estimation and the theoretical approximation (\ref{theoric}), for the AR(1)-input with $b=3$ and $n=10^4$}
\end{table}

\newpage
\begin{center}\large\scshape Appendix. $\tau$-weak dependence and proofs\end{center}
\appendix
\begin{itemize}\item[A.1.] {\it \large Brief interlude into $\tau$-weak dependence}\end{itemize}

Let $(\Omega, \mc{A}, \PP)$ be a probability space, and $\mc{M}$ a $\sigma$-algebra of $\mc{A}$. Let $(E, \delta)$ be a Polish space endowed with its metric. For any $E$-valued random variable $X$, $\mathbb{L}^p$-integrable ({\it i.e.} $X$ satisfies $\|\delta(X,0)\|_p:=\left(\int \delta^p(x,0)\PP(dx)\right)^{1/p}<\infty$),   \citet{Dedecker2004a} defined the \textbf{coefficient $\tau_p$} as:
\begin{eqnarray}\label{tau1}
\tau_p(\mc{M}, X):=\|\sup\left\{\EE\left[h(X)|\mc{M}\right]-\EE\left[h(X)\right]: h\in\Lambda(E,\delta)\right\}\|_p
\end{eqnarray}
where $\Lambda(E,\delta)$ denotes the class of all Lipschitz functions $h:E\lfled\RR$ such that $$\Lip(h):=\sup_{x\­y}\frac{|h(x)-h(y)|}{\delta(x,y)}\le 1.$$

Let $\mc{X}=(X_{n,i})_{1\²i\²n, n\in\NN}$ be a triangular array of\, $\mb{L}^p$-integrable $E$-valued random variables, and $(\mc{M}_{i})_{i\in\ZZ}$  be a sequence of $\sigma$-algebras of $\mc{A}$. \\
 Then, for any $n\in\NN$, we define the coefficient: 
\begin{align}\label{tau2}
\tau_{p,n}(k):=\sup_{l\³1} l^{-1}\sup\{ \tau_p(\mc{M}_{i},(X_{n,j_1}, \ldots, X_{n,j_l})): i+k\² j_1< \cdots < j_l \² n\},
\end{align}
where we  consider the distance 
\begin{align}\label{metric}
\delta_l(x,y)=\sum_{i=1}^l \delta(x_i, y_i)
\end{align}
on $E^l$.  Moreover, we say that $\mc{X}$ is $\tau_p$-weakly dependent if 
\begin{align}\label{tau3}
\lim_{k\lfled \infty} \limsup_{n\lfled \infty} \tau_{p,n}(k)=0. 
\end{align}

\begin{rem}\rm Recall that $\mc{X}=(X_{n,i})_{1\²i\²n, n\in\NN}$ is constructed from a random process $\mathbb{X}=(X_i)_{i\in\ZZ}$. Therefore the dependence properties of $\mc{X}$'s  are inherited from those of  $\mathbb{X}$. Even more so, if $\mathbb{X}$ is $\tau_p$-weakly dependent (in the usual sense defined by \citet{Dedecker2004a} for random processes), then $\mc{X}$ is $\tau_p$ - weakly dependent with $\tau_{p,n}(\cdot)=L_n\cdot \tau_p(\cdot)$, for some positive constant $L_n$ (which is written in function of  $\mc{X}$'s normalization constants). For instance, if we consider the normalizations in eqn. (\ref{N1})-(\ref{4}) we obtain that $\tau_{p,n}(\cdot)\²a^{-1}_n\cdot \tau_p(\cdot)$.

In this sense, if we want to study $\mc{X}$'s $\tau$-dependence properties, suffice it to take into account $\mathbb{X}$'s $\tau$-dependence properties. 
\end{rem}

We make use of the previous remark to mention the following examples of $\tau_p$-weakly dependent processes without considering the normalizations. 

\begin{exa}[Causal Bernoulli shifts]\rm Let $(\xi_i)_{i\in\ZZ} $ be a sequence of i.i.d.r.v's. (independent and identically distributed random variables) with values in a measurable space $D$. Assume that there exists a function $H: D^{\NN}\lfled \RR$, such that $H(\xi_0,\xi_{-1},\ldots )$ is defined almost surely. Then the stationary sequence $(X_i)_{i\³0}$ defined by $X_i=H(\xi_i,\xi_{i-1},\ldots)$ is called a causal Bernoulli shifts. 

Let $(\xi'_i)_{i\in\ZZ}$ be an independent copy  of the  i.i.d. sequence $(\xi_i)_{i\in\ZZ}$.  Consider a decreasing sequence  $(\Delta_p(i))_{i\³0}$ such that 
\begin{align}\label{N0}
\|X_i-X'_i\|_p\²\Delta_p(i)
\end{align}
for some $p\in[1,\infty]$, where $X_i=H(\xi_i,\xi_{i-1}, \ldots)$ and $X'_i=H(\xi_i, \ldots, \xi_1, \xi'_0,\xi'_{-1},\ldots)$. Then, if $\mc{M}_i=\sigma(X_j: j\² i)$, the coefficient $\tau_p(k)$ of $(X_i)_{i\in\ZZ}$ is bounded above by $\Delta_p(k)$, for all $k\in\NN$.
\end{exa}
\begin{app}[Causal linear processes]  Let $D=\RR$ and
\begin{align}\label{CLP}
X_i=\sum_{j=0}^\infty b_j \xi_{i-j}.
\end{align} Here we  set $\Delta'_p(i)=2\|\xi_0\|_p\sum_{j\³i} |b_j|\³ \Delta_p(i)$ in case $\|\xi_0\|_p<\infty$.\\
 The model (\ref{AR}) can be write as (\ref{CLP}) with  $b_j=b^{-j-1}$ for some integer $b\³2$. Since $\xi_0$ is  uniformly distributed on $U(b)=\{0, \ldots, b-1\}$ in this case, then $X_0$ is uniformly distributed over $[0,1]$ and $\Delta_\infty (i)\²b^{-i}$.  
\end{app}
\begin{exa}[Markov models]\rm Let $G:(\RR^l,\mc{B}(\RR^l)) \times (D, \mc{D})\lfled (\RR, \mc{B}(\RR))$ be a measurable function and let $(X_i)_{i\³ 1-l}$ be a sequence of random variables with values in $\RR$ such that 
\begin{align}\label{MC}
X_i=G(X_{i-1}, X_{i-2},\ldots, X_{i-l}; \xi_i) ,\qquad \forall i\³1,
\end{align} 
for some sequence $(\xi_i)_{i\in\NN}$ of i.i.d.r.v's. with values in a measurable space $D$ and independent of $(X_0, \ldots, X_{1-l})$. Then the random variables $Y_i=(X_i, X_{i-1}, \ldots, X_{i-l+1})$ defines a Markov chain such that $Y_i=F(Y_{i-1};\xi_i)$ with
\begin{align}\label{MC2}
F(x_{l},\ldots, x_1; \xi):=(G(x_{l}, \ldots, x_1;\xi), x_{l}, x_{l-1}, \ldots, x_2).
\end{align}
Assume that $Y_0=(X_0, \ldots, X_{1-l})$ is a stationary solution of (\ref{MC}). Let $Y'_0=(X'_0,\ldots, X'_{1-l})$ be independent of  $(Y_0, (\xi_i)_{i\in\NN})$ and distributed as $Y_0$. Then setting 
\begin{align}\label{MC3}
X'_i=G(X'_{i-1},\ldots, X'_{i-l}; \xi_i),
\end{align}  $X'_i$ is distributed as $X_i$ and it is independent of $\mc{M}_0=\sigma(X_0, \ldots, X_{1-l})$, for all $i\in\NN$. As in the previous example, let $(\Delta_{p}(i))_{i\³0}$ be a non increasing sequence such that (\ref{N0}) holds, with $X_i$ and $X'_i$ defined in (\ref{MC}) and (\ref{MC3}), respectively. Hence one can apply the result of Lemma 3 in \citep{Dedecker2004a}, and we obtain that $\tau_p(k)\² \Delta_p(k)$. 

In particular if $G$ is such that 
\begin{align}\label{Stat}
\|G(x; \xi_1)-G(y;\xi_1)\|_p\² \sum_{i=1}^l a_i | x_i-y_i|\text{, \quad with } \sum_{i=1}^l a_i <1,
\end{align}
then $\Delta_p(k)\² C a^k$ for some $a\in[0,1)$ and some $C>0$. (see \citep{Dedecker2007}, page 34).
\end{exa}

\begin{app}[Contracting Markov chain]  Let $X_i=G(X_{i-1},\xi_i)$ be a Markov chain such that $G: (\RR, \mc{B}(\RR))\times (D, \mc{D})\lfled (\RR, \mc{B}(\RR))$ is a measurable function and 
\begin{align}\label{stat2}
A=\|G(0;\xi_1)\|_p<\infty \text{ and } \|G(x; \xi_1)-G(y;\xi_1) \|_p\² a|x-y|,
\end{align} 
for some $a\in (0,1)$ and some $p\in [1,\infty]$. Then, $(X_i)_{i\in\NN}$ has a stationary solution with $p$-th order finite moment as this is proved on page 35 of \citet{Dedecker2007}. Moreover under this condition: $\Delta_p(k)=\|X'_0-X_0\|_p\cdot  a^k$. 
\end{app}

\begin{rem}\rm In particular if $G(u;\xi)=A(u)+B(u)\xi$ for suitable Lipschitz functions $A(u)$ and $ B(u)$ with $u\in\RR$, then the corresponding iterative model (ARCH-type process) $X_i=G(X_{i-1};\xi_i)$ satisfies (\ref{stat2}) with $a=\Lip(A)+\|\xi_1\|_p \Lip(B)<1$.
\end{rem}

\begin{rem}\rm The stationary iterative models $X_i=G(X_{i-1},\xi_i)$ are causal Bernoulli shifts if the condition (\ref{stat2}) holds; this is proved  in Proposition 3.2 in \citep{Dedecker2007}.
\end{rem}

\begin{app}[Nonlinear AR($l$)-models] Let $l\ge1$ and $(X_i)_i $ be the stationary solution of some equation  $$X_i=R(X_{i-1},\ldots, X_{i-l})+\xi_i$$ 
for some measurable function  $R:\RR^l\to\RR$. The process  $(X_i)_i $  is then called  a stationary real nonlinear autoregressive model of order $l$. If $\|\xi_1\|_p<\infty$ and $$|R(u_1,\ldots, u_l)-R(v_1,\ldots, v_l)|\²\sum_{i=1}^l a_i |u_i-v_i|, \mbox{ for  }a_1, \ldots, a_l\³0\mbox{ with } \sum_{i=1}^l a_i<1,$$ and for all $(u_1,\ldots,u_l),(v_1,\ldots,v_l)\in\RR^l$, then the function $G:\RR^{l+1}\to \RR$ defined by  $G(u; \xi)=R(u)+\xi$ satisfies Condition (\ref{Stat}) and therefore the sequence $(\Delta_p(k))_k$ admits an  exponential decay rate. 
\end{app}
\begin{itemize}\item[A.2.] {\it \large Coupling}\end{itemize}
Let  $d_\cup: E_\cup \times E_\cup \lfled [0,\infty)$ be a pseudometric defined as follows:
let $x_k=(x_k^{(1)}, x_k^{(2)}, \ldots , x_{k}^{(r_k)}) \in E_\cup$ for $k=1,2$ and we denote 
\[x_{_{1,2}}=\left\{\begin{array}{ll}
x_1, & \mbox{if $r_1\³ r_2$} ,\\ 
(x_1^{(1)},\ldots, x_{1}^{(r_1)}, x_2^{(r_1+1)}, \ldots, x_2^{(r_2)}), & \mbox{if  $r_1< r_2$}. 
\end{array}\right.\]
Similarly is denoted $x_{_{2,1}}$.
Then,
\begin{align}\label{semimetric}
d_\cup(x_1, x_2):=d_{r_1\vee r_2} (x_{_{1,2}} , x_{_{2,1}}),
\end{align}
where $d_r(\cdot , \cdot)$ is defined in (\ref{metric}).

\begin{lem}[Coupling: the even and odd blocks] 
Suppose that the random variables $(X_{n,i})_{1\²i\²n}$ are such that (B.2) holds. We consider together even and odd block sizes, by using $k=0$ or $1$ according to the parity. Assume that for each $j\in\{2,\ldots, [m_n/2]\}$ there is a random variable $U_{k,j}$ uniformly distributed on $[0,1]$ and independent of  $\mc{M}_{n,j-1}^{k}=\sigma(Y_{n,2-k}, \ldots, Y_{n,2(j-1)-k)})$ and $\sigma(Y_{n,2j-k})$. Then there exists a random block $\dot{Y}_{n,2j-k}$ measurable with respect to $\mc{M}_{n,j-1}^k\vee \sigma(Y_{n,2j-k})\vee \sigma(U_{k,j})$, independent of $\mc{M}_{n,j-1}^k$ and distributed as $Y_{n,2j-k}$ such that 
\begin{align}\label{coutau1}
\left\|\left.\EE\left(d_\cup (Y_{n,2j-k},\dot{Y}_{n,2j-k})\right)\right|\mc{M}_{n,j-1}^k\right\|_1\²r_n \tau_{1,n}(r_n).
\end{align}
In particular, if we set  $\mc{M}_{n,j-1}^k=\sigma(\dot{Y}_{n,2-k},\ldots, \dot{Y}_{n,2(j-1)-k})$, then the blocks $(\dot{Y}_{n,2j-k})_{1\²j\²[m_n/2]}$ are independent. \\\end{lem}
\noindent
{\it \textbf{Proof:}}
We set here  $k=0$ (for even block sizes) since the steps are similar if $k=1$.\\ 
Let $Y_{n,2j}=(X_{n,i})_{i\in A_j}$ be a random block, where$A_j=\{(2j-1)r_n +1, \ldots, 2jr_n \}$. Then, from \citet{Dedecker2007}'s Lemma 5.3, there exists a random random block $\dot{Y}_{n,2j}\in E^{r_n}$ measurable with respect to $\mc{M}_{n,j-1}^0 \vee \sigma(Y_{n,2j})\vee \sigma(U_{0,j})$ independent of $\mc{M}_{n,j-1}^0$  and distributed as $Y_{n,2j}$  such that 
\begin{align}
\tau_p(\mc{M}^0_{n,j-1}, Y_{n,2j})=\left\| \EE ( d_{r_n}(Y_{n,2j}, \dot{Y}_{n,2j})\left| \mc{M}^0_{n,j-1})\right.\right\|_p.
\end{align}
Therefore, 
\begin{multline}
\left\|\EE( d_{\cup}\left( Y_{n,2j}, \dot{Y}_{n,2j}\right) \left| \mc{M}^0_{n,j-1})\right.\right\|_1
=\left\|\EE( d_{r_n}\left( Y_{n,2j}, \dot{Y}_{n,2j}\right) \left| \mc{M}^0_{n,j-1})\right.\right\|_1
\\
=\tau_1(\mc{M}^0_{n,j-1}, Y_{n,2j})
\le  r_n \tau_{1,n}(r_n).\ \  \square
\end{multline}

\begin{lem}[Coupling: the sub-blocks] 
Suppose that the random variables $(X_{n,i})_{1\²i\²n}$ are such (B.2) holds. Moreover, if that for each $j\in\{2,\ldots, [m_n/2]\}$ there is a random variable $U_{k,j}$ uniformly distributed on $[0,1]$ and independent of the $\sigma-$algebras $\mc{M}_{n,j-1}=\sigma(Y_{n,1}^{[r_n-l_n]}, \ldots, Y_{n,j-1}^{[r_n-l_n]})$ and   $\sigma(Y_{n,j}^{[r_n-l_n]})$. Then there exists a random block $\dot{Y}_{n,j}^{[r_n-l_n]}$, measurable with respect to $\mc{M}_{n,j-1}\vee \sigma(Y_{n,j}^{[r_n-l_n]})\vee \sigma(U_j)$, independent of $\mc{M}_{n,j-1}$ and distributed as $Y_{n,j}^{[r_n-l_n]}$ such that  
\begin{align}\label{coutau2}
\left\|\left.\EE\left(d_\cup(Y_{n,j}^{[r_n-l_n]},\dot{Y}_{n,j}^{[r_n-l_n]})\right)\right|\mc{M}_{n,j-1}\right\|_1\²r_n\tau_{1,n}(l_n).
\end{align} 
Moreover, if we set $\mc{M}_{n,j-1}=\sigma(\dot{Y}_{n,1}^{[r_n-l_n]},\ldots, \dot{Y}_{n,j-1}^{[r_n-l_n]})$ then the blocks $(\dot{Y}_{n,j}^{[r_n-l_n]})_{1\²j\²m_n}$ are independent. 
\end{lem}\noindent
{\it Proof.}
The same argument previous proof. However note that the sub-blocks $(Y_{n,j}^{[r_n-l_n]})_{1\²j\²m_n}$ are separated by $l_n$ variables.
\hfill$\square$\vskip2mm\hfill
\begin{itemize}\item[A.3.] {\it \large Proof of Theorem 1}\end{itemize}
Let $(Y_{n,j})_{1\²j\²m_n}$ be the blocks built from $(X_{n,i})_{1\²i\²n}$. For $k\in\{0,1\}$, we consider the independent blocks $(\dot{Y}_{n,2j-k})_{1\²j\²[m_n/2]}$ coupled to the original blocks $(Y_{n,2j-k})_{1\²j\²[m_n/2]}$, from Lemma 2. Therefore, if we define $\Delta_{n,j}^*:=f(\dot{Y}_{n,j})-f(\dot{Y}_{n,j}^{(r_n-l_n)})$, for $j=1,\ldots, m_n$, we have that $\Delta_{n,j}^*(f)\overset{\mc{D}}{=}\Delta_{n,j}(f)\overset{\mc{D}}{=}\Delta_n(f)$, for each $j$, where $\Delta_{n,j}(f):=f(Y_{n,j})-f(Y_{n,j}^{(r_n-l_n)})$ and $\Delta_n(f)$ is defined in (\ref{Delta}). Now, if we consider the assumption (C.1), we can apply  \citet{Petrov1975}'s Theorem 1 (Section IX.1) to the i.i.d.r.v's $X_{n,j}:=(nv_n)^{-1/2}\Delta_{n,j}^*(f)$, so
\begin{align}\label{Delta2}
D\dot{Z}_n^{(k)}(f):=\dfrac{1}{\sqrt{n v_n}}\sum_{j=1}^{[m_n/2]} \left(\Delta^*_{n,2j-k}(f)-\EE\Delta^*_{n,2j-k}(f)\right)=o_P(1)
\end{align}
for $k=0,1$. In consequence,
\begin{align}\label{Delta3}
DZ_n(f):=\dfrac{1}{\sqrt{n v_n}}\sum_{j=1}^{m_n} \left(\Delta_{n,j}(f)-\EE\Delta_{n,j}(f)\right)=o_P(1)
\end{align}

On the other hand, by Lemma 3, we have that 
\begin{align}\label{B1}
BZ_n(f):=\dfrac{1}{\sqrt{n v_n}}\sum_{j=1}^{m_n} \left(f(Y_{n,j}^{(r_n-l_n)})-\EE f(Y_{n,j}^{(r_n-l_n)})\right)
\end{align}
converge weakly in fidis if, and only if 
\begin{align}\label{B2}
B\dot{Z}_n(f):=\dfrac{1}{\sqrt{n v_n}}\sum_{j=1}^{m_n} \left(f(\dot{Y}_{n,j}^{(r_n-l_n)})-\EE f(\dot{Y}_{n,j}^{(r_n-l_n)})\right)
\end{align}
converge weakly in fidis (and in this case the limit distributions are the same).  The latter holds because $B\dot{Z}_n(f)=\dot{Z}_n(f)-D\dot{Z}_n(f)$ and from the assumptions (C.2) and (C.3), where 
$$\dot{Z}_n(f):= \dfrac{1}{\sqrt{nv_n}}\sum_{j=1}^{m_n}\left(f(\dot{Y}_{n,j})-\EE f(\dot{Y}_{n,j})\right).$$

Finally, as $Z_n(f)=BZ_n(f)+DZ_n(f)$  $\forall f\in\mc{F}$, we get the result.
\hfill$\square$\vskip2mm\hfill
\begin{itemize}\item[A.4.] {\it \large Proof of Proposition 1}\end{itemize}
Suffices to prove the following multidimensional version of \citet{Segers2003}'s condition (6):
\begin{align}\label{Segerscondition}
\lim_{l\to\infty}\limsup_{n\to\infty} \PP\left\{Y_n^{[l+1: r_n]}\­0 | X_{n,1}\­0 \right\}=0, 
\end{align}
since the rest of the proof follows the same steps of the proof of \citet{Drees2010}'s Lemma 5.2.

Indeed, let $h(\cdot)=\1\{ \cdot \­ 0\}$ be a function defined on $E^{r_n -l }$. Consider a increasing sequence of functions $h_k(\cdot):E^{r_n-l}\lfled [0,1]$ which approximate to $h$, and such that $\Lip(h_k)=v_k^{-\alpha}$ for some $\alpha>0$. Of course, we set $k=k(n) \ll  n$. Then, 
\begin{multline*}
\limsup_{n\to\infty}\PP\{ Y_n^{(l+1: r_n)}\­0| X_{n,1}\­0\}= \limsup_{n\to\infty}v_n^{-1} \PP\{Y_n^{(l+1: r_n)}\­0,X_{n,1}\­0\}
\\ =\limsup_{n\to\infty}v_n^{-1}\int_{\{X_{n,1}\­0\}}\PP\{Y_n^{(l+1: r_n)}\­0|\sigma(X_{n,1})\}d\PP
\\=\limsup_{n\to\infty} v_n^{-1}\int_{\{X_{n,1}\­0\}}\EE[\mb{1}\{Y_n^{(l+1: r_n)}\­0\}|\sigma(X_{n,1})]d\PP
\\ =\limsup_{k\to\infty}\dfrac{v_n^{-1}}{v_k^\alpha}\int_{\{X_{n,1}\­0\}}\EE\left[\dfrac{h_k(Y_n^{[l+1:r_n]})}{v_k^{-\alpha}}|\sigma(X_{n,1})\right]d\PP
\\ \² \limsup_{k\to\infty}\dfrac{1}{v_k^{\alpha+1/p}}\left\| \EE\left[\dfrac{h_k(Y_n^{[l+1:r_n]})}{v_k^{-\alpha}}|\sigma(X_{n,1})\right]\right\|_p
\\ \² \limsup_{k\to\infty}\left[\dfrac{(r_n-l)}{v_k^{\alpha+1/p}}\tau_{n,p}(l)+ (r_n-l) v_n^{1/q}\right]\² \limsup_{k\to\infty}\dfrac{r_n\tau_{n,p}(l)}{v_k^{\alpha+1/p}}.
\end{multline*}
Finally, taking $l\to\infty$ we have the limit (\ref{Segerscondition}) proven.
\hfill$\square$\vskip2mm\hfill
\begin{itemize}\item[A.5.] {\it \large Proof of Theorem 2}\end{itemize}
Note that $Z_{n}$ is asymptotically tight iff $\dot{Z}_n^{(k)}$ defined by:
\begin{align}\label{even-odd1}
Z_n^{(k)}(f):=\dfrac{1}{\sqrt{nv_n}} \sum_{j=1}^{[m_n/2]} \left(f(Y_{n,2j-k})-\EE f(Y_{n,2j-k})\right)
\end{align} 
is asymptotically tight for each $k\in\{0,1\}$.  On the other hand, for each $k\in\{0,1\}$ we use Lemma 2 together with (B.2) and (D.1) conditions to build independent blocks $(\dot{Y}_{n,2j-k})_{1\²j\²[m_n/2]}$ coupled to the original blocks $(Y_{n,2j-k})_{1\²j\²[m_n/2]}$. In this manner we have that $Z_{n}^{(k)}$ is asymptotically tight iff 
\begin{align}\label{even-odd2}
\dot{Z}_n^{(k)}(f):=\dfrac{1}{\sqrt{nv_n}} \sum_{j=1}^{[m_n/2]} \left(f(\dot{Y}_{n,2j-k})-\EE f(\dot{Y}_{n,2j-k})\right)
\end{align} 
is asymptotically tight, for each $k\in\{0,1\}$. The latter is true due to Theorem 2.11.9 in \citep{VdV} by setting $Z_{nj}(f)=f(Y_{n,j})$ and $[m_n/2]$ instead of $m_n$. 

For the remaining assertion we use Theorem 1.
\hfill$\square$\vskip2mm\hfill
\begin{itemize}\item[A.6.] {\it \large Proof of Theorem 3}\end{itemize}
Consider (T.5). Note that  from the triangle inequality $Z_n$ is asymptotically equicontinuous if $\dot{Z}^{(k)}_n$ from eqn. (\ref{even-odd1}) is asymptotically equicontinuous for each $k\in\{0,1\}$. Now, again we use Lemma 2 together with (B.2) and (D.1) conditions as in the previous proof to prove that $Z_{n}^{(k)}$ is asymptotically equicontinuous iff  $\dot{Z}_n^{(k)}$ is asymptotically equicontinuous for each $k\in\{0,1\}$. However, in this case $\dot{Z}_n^{(k)}$ is  asymptotically equicontinuous from Theorem 2.11.1 in \citep{VdV}.
\\
The remaining steps are the same of Theorem 2.10's proof in \citep{Drees2010}. 
\hfill$\square$\vskip2mm\hfill
\begin{itemize}\item[A.7.] {\it \large Proof of Proposition 2}\end{itemize}
 The steps are the same that in the proof of Theorem 2.1 in \citep{Drees2011}, but replacing the assumptions (C1) and (C2) of his paper by our assumptions (B.2) and (B.1), respectively.
 
\hfill$\square$\vskip2mm\hfill
\begin{itemize}\item[A.8.] {\it \large Proof of Corollary 3}\end{itemize}
Suffices to replace the assumptions (C1) and (C2) in the proof of \citet{Drees2011}'s Corollary 2.3  by our assumptions (B.2) and (B.1) respectively.
\hfill$\square$\vskip2mm\hfill
\begin{itemize}\item[A.9.] {\it \large Proof of the expression (\ref{dgn})}\end{itemize}
If $(X_i)_{i\³0}$ is the AR(1)-process (\ref{AR}), note that for each $i\in\NN$
\begin{align}
X_i= \dfrac{X_0}{b^i}+\sum_{s=1}^i \dfrac{\xi_s}{b^{i-s+1}}.
\end{align} 
Therefore, if $n$ is sufficiently large such that $b^d v_n<1$, then for $x\in [0,1]^d$:
\begin{multline*}
\PP\{X_{n,1}> x, X_{n,1}\­0\}=\PP\{ X_i>a_n x_i + u_n, \text{ for some } i=1,\ldots, d\}
\\ =\PP\left\{X_0>b^i(a_n x_i+u_n)-\sum_{s=1}^i b^{s-1} \xi_s,  \text{ for some  } i=1, \ldots,d \right\}
\end{multline*}

\begin{multline*}
=\PP\left\{X_0>\min_{i=1, \ldots, d}\left\{b^i(a_n x_i+u_n)-\sum_{s=1}^i b^{s-1} \xi_s\right\}\right\}
\\=\sum_{j_1, \ldots, j_d \in U(b)}\PP\left\{X_0>\min_{i=1, \ldots, d}\left\{b^i(a_n x_i+u_n)-\sum_{s=1}^i b^{s-1} j_s\right\}, (\xi_1, \ldots, \xi_d)=(j_1, \ldots, j_d)\right\}
\\=\dfrac{1}{b^d}\sum_{j_1, \ldots, j_d \in U(b)}\PP\left\{X_0>\min_{i=1, \ldots, d}\left\{b^i(a_n x_i+u_n)-\sum_{s=1}^i b^{s-1} j_s\right\}\right\}
\\=\dfrac{1}{b^d}\sum_{j_1, \ldots, j_d \in U(b)}\left(\max_{i=1, \ldots, d}\left\{1-b^i+\sum_{s=1}^i b^{s-1}j_s+b^i v_n (1-x_i)\right\}_+ \land 1\right)
\\=\dfrac{1}{b^d}\max_{i=1, \ldots, d}\left\{b^i v_n (1-x_i)\right\},
\end{multline*}
since $\mu_b(j_1,\ldots, j_d; i):=1-b^i+\sum_{s=1}^i b^{s-1}j_s\² -1$ for all $(j_1,\ldots, j_d)\in U^d(b)\setminus \{(b-1, \ldots, b-1)\}$ and $\mu_b(b-1, b-1, \ldots, b-1)=0$.
Thus, 
\begin{align}
\PP\{X_{n,1}> x| X_{n,1}\­0\} \underset{_{n\to \infty}}{\lfled}  \max_{i=1,\ldots, d}\{ b^{i-d}(1-x_i)\}.
\end{align}
\hfill$\square$\vskip2mm\hfill
\begin{itemize}\item[A.10.] {\it \large Proof of the expression (\ref{covar}})\end{itemize}
Let $x, y\in [0,1]^d$. Then as before for $i\³1$, if $n$ is sufficiently large such that $b^{i+d} v_n<1$, then we have:
\begin{multline*}
\PP\{X_{n,1}>x, X_{n,i+1}>y\}
\\=\PP\left\{ X_k> a_n x_k + u_n, X_{i+l}> a_n y_l + u_n, \text{ for some } (k,j)\in\{1, \ldots, d\}^2\right\}
\\=\PP\left\{ X_0> \min_{k=1,\ldots, d}\left\{ b^k(a_n x_k + u_n)-\sum_{s=1}^k b^{s-1} \xi_s\right\},\right.
\\ \left. X_i> \min_{l=1, \ldots, d}\left\{ b^l(a_n y_l + u_n) -\sum_{s=1}^l b^{s-1} \xi_{s+i}\right\}\right\}
\\=\sum_{\overset{j_1,\ldots, j_d\in U(b)}{ j_{i+1},\ldots, j_{i+d}\in U(b)}}\PP\left\{ X_0> \min_{k=1,\ldots, d}\left\{ b^k(a_n x_k + u_n)-\sum_{s=1}^k b^{s-1} j_s\right\},\right.
\\ X_i> \min_{l=1, \ldots, d}\left\{ b^l(a_n y_l + u_n) -\sum_{s=1}^l b^{s-1} j_{s+i}\right\}, 
\\ \left.(\xi_1, \ldots, \xi_d, \xi_{i+1},\ldots, \xi_{i+d})=(j_1, \ldots, j_d, j_{i+1}, \ldots, j_{i+d})\right\}
\end{multline*}
\begin{multline*}
=\dfrac{1}{b^d}\sum_{\overset{j_1,\ldots, j_d\in U(b)}{ j_{i+1},\ldots, j_{i+d}\in U(b)}}\PP\left\{ X_0> \min_{k=1,\ldots, d}\left\{ b^k(a_n x_k + u_n)-\sum_{s=1}^k b^{s-1} j_s\right\},\right.
\\ \left. X_i> \min_{l=1, \ldots, d}\left\{ b^l(a_n y_l + u_n) -\sum_{s=1}^l b^{s-1} j_{s+i}\right\}, (\xi_1, \ldots, \xi_d)=(j_1, \ldots, j_d)\right\}
\\=\dfrac{1}{b^d}\sum_{\overset{j_1,\ldots, j_d\in U(b)}{ j_{i+1},\ldots, j_{i+d}\in U(b)}}\PP\left\{ X_0>1- \max_{k=1,\ldots, d}\left\{ \mu_b(j_1, \ldots, j_d; k)+ b^kv_n(1-x_k)\right\}_+\land 1,\right.
\\ \left. X_i> 1- \max_{l=1, \ldots, d}\left\{ \mu_b(j_{i+1,\ldots, j_{i+d}})+ b^lv_n(1-y_l)\right\}_+\land 1, (\xi_1, \ldots, \xi_d)=(j_1, \ldots, j_d)\right\}%
\\=\dfrac{1}{b^d}\PP\left\{ X_0>1- \max_{k=1,\ldots, d}\left\{b^k v_n (1-x_k)\right\},\right.
\\ \left. X_i> 1- \max_{l=1, \ldots, d}\left\{ b^l v_n (1-y_l)\right\}, \xi_1=\ldots=\xi_d=b-1\right\}
\end{multline*}
since $\mu_b(j_1,\ldots, j_d; i):=1-b^i+\sum_{s=1}^i b^{s-1}j_s\² -1$ for all $(j_1,\ldots, j_d)\in U^d(b)\setminus \{(b-1, \ldots, b-1)\}$ and $\mu_b(b-1, b-1, \ldots, b-1)=0$.

Moreover, note that if $i>d$ 
\begin{multline*}
\PP\{X_{n,1}>x, X_{n,i+1}>y\}=\dfrac{1}{b^d}\PP\left\{ X_0>1- \max_{k=1,\ldots, d}\left\{b^k v_n (1-x_k)\right\},\right.
\\ \left. X_i> 1- \max_{l=1, \ldots, d}\left\{ b^l v_n (1-y_l)\right\}, \xi_1=\ldots=\xi_d=b-1\right\}
\\=\dfrac{1}{b^d}\PP\left\{ X_0>1- \max_{k=1,\ldots, d}\left\{b^k v_n (1-x_k)\right\},\right.
\\ \left. X_0> b^i-b^i\max_{l=1, \ldots, d}\left\{ b^l v_n (1-y_l)\right\}+1-b^d-\sum_{s=d+1}^i b^{s-1} \xi_s\right\}
\\ =\dfrac{1}{b^d}\sum_{j_{d+1, \ldots, j_i}\in U(b)}\PP\left\{ X_0>1- \max_{k=1,\ldots, d}\left\{b^k v_n (1-x_k)\right\},\right.
\\ \left. X_0> b^i-b^i\max_{l=1, \ldots, d}\left\{ b^l v_n (1-y_l)\right\}+1-b^d-\sum_{s=d+1}^i b^{s-1} j_s, (\xi_{d+1},\ldots, \xi_i)=(j_{d+1}, \ldots, j_i)\right\}
\\=\dfrac{1}{b^i}\sum_{j_{d+1, \ldots, j_i}\in U(b)}\PP\left\{ X_0>1- \max_{k=1,\ldots, d}\left\{b^k v_n (1-x_k)\right\},\right.
\\ \left. X_0> b^i-b^i\max_{l=1, \ldots, d}\left\{ b^l v_n (1-y_l)\right\}+1-b^d-\sum_{s=d+1}^i b^{s-1} j_s\right\}
\end{multline*}
\begin{multline*}
=\dfrac{1}{b^i}\sum_{j_{d+1, \ldots, j_i}\in U(b)}\min\left(\max_{k=1,\ldots, d}\left\{b^k v_n (1-x_k)\right\},\right.
\\ \left.\max_{l=1, \ldots, d}\left\{ b^d+\sum_{s=d+1}^i b^{s-1} j_s + b^{l+i} v_n (1-y_l)-b^i\right\}_+, 1\right)\\ \\ =\dfrac{v_n}{b^i}\min\left(\max_{k=1,\ldots, d}\left\{b^k (1-x_k)\right\}, \max_{l=1, \ldots, d}\left\{b^{l+i}(1-y_l)\right\}\right)= v_n H_{b,i}(x,y)
\end{multline*}
 Similarly for $1\²i< d$, we obtain that
 \begin{multline*}
\PP\{X_{n,1}>x, X_{n,i+1}>y\}
\\=\dfrac{v_n}{b^i}\min\left( \max_{k=1\ldots, i} \{ b^k(1-x_k)\}, \max_{k=i+1, \ldots, d}\{b^k\min(1-x_k, 1-y_k)\}, \max_{k=d-i, \ldots, d} \{ b^{k+i} (1-y_k)\}\right)
\\=v_n H_{b,i}(x,y)
\end{multline*}
From Lemma 5.2 - (iii) in \citep{Drees2010}, $\EE |f(Y_n)|=o(\sqrt{nv_n})$. Thus, for $n$ sufficiently large: 
\begin{multline*}
\dfrac{\Cov \left(f_x(Y_n), f_y (Y_n)\right)}{r_n v_n} \sim \PP\{X_{n,1}>x, X_{n,1}>y\}
\\ + \sum_{i=1}^{r_n-1}\left(1-\frac{i}{r_n}\right) \left(\PP\{X_{n,1}>x, X_{n,i+1}>y\} + \PP\{X_{n,1}>y, X_{n,i+1}>x\}\right) 
\\ \underset{_{n\to \infty}}{\lfled} \min\left( \max_{k=1,\ldots, d}\{ b^k (1-x_k)\}, \max_{k=1, \ldots, d}\{ b^k(1-y_k)\}\right) 
+\sum_{i=1}^\infty \left(H_{b,i}(x,y) + H_{b,i}(y,x)\right).
\end{multline*}
\hfill$\square$\vskip2mm\hfill
\begin{itemize}\item[A.11.] {\it \large Proof of the expression (\ref{theoric}})\end{itemize}
The proof is similar to the proof of the expression (\ref{dgn}). Indeed, 
\begin{multline*}
\PP\left\{ \max_{1\²i\²r_n} X_i> 1-v_n t \right\}=\PP\left\{X_i> 1-v_n t, \text{ for some } i=1, \ldots, r_n \right\}
\\ = \PP\left\{ X_0> b^i(1-v_n t)-\sum_{s=1}^i b^{s-1}\xi_s, \text{ for some  } i=1, \ldots, r_n\right\}
\\ = \PP\left\{ X_0> \min_{1\²i\²r_n}\left\{b^i(1-v_n t)-\sum_{s=1}^i b^{s-1}\xi_s\right\}\right\}
\end{multline*}
\begin{multline*}
=\sum_{j_1, \ldots, j_{r_n}\in U(b)}\PP\left\{ X_0> \min_{1\²i\²r_n}\left\{b^i(1-v_n t)-\sum_{s=1}^i b^{s-1}\xi_s\right\}, (\xi_1, \ldots, \xi_{r_n})=(j_1, \ldots, j_{r_n})\right\}
\\= \dfrac{1}{b^{r_n}}\sum_{j_1, \ldots, j_{r_n}\in U(b)}\PP\left\{ X_0> \min_{1\²i\²r_n}\left\{b^i(1-v_n t)-\sum_{s=1}^i b^{s-1}j_s\right\}\right\}
\\= \dfrac{1}{b^{r_n}}\sum_{j_1, \ldots, j_{r_n}\in U(b)}\min\left(\max_{1\²i\²r_n}\left\{1+\sum_{s=1}^i b^{s-1}j_s - b^i(1-v_nt)\right\}_+, 1\right).
\end{multline*}
\hfill$\square$\vskip2mm\hfill
\paragraph{\bf Aknowledgements.} 
Warm thanks are due to the very constructive and friendly help of Olivier Wintenberger in the redaction and the finalization of this paper.\\
 Very specials thanks are also due to an anonymous referee who pointed clearly the weaknesses of a previous version of this work and helped us to make it more adequate for publication. 



\end{document}